\documentclass{amsart}
\usepackage{amssymb} 
\usepackage{amsmath} 
\usepackage{amscd}
\usepackage{amsbsy}
\usepackage{comment}
\usepackage[matrix,arrow]{xy}
\usepackage{hyperref}

\DeclareMathOperator{\Norm}{Norm}

\DeclareMathOperator{\Res}{Res}
\DeclareMathOperator{\rank}{rank}

\DeclareMathOperator{\Gal}{Gal}

\newcommand{\Q}{{\mathbb Q}}
\newcommand{\Z}{{\mathbb Z}}
\newcommand{\F}{{\mathbb F}}

\newcommand{\cN}{\mathcal{N}}

\newcommand{\OO}{{\mathcal O}}

\newcommand{\fp}{\mathfrak{p}}
\newcommand{\fq}{\mathfrak{q}}

\newcommand{\bs}{\mathbf{s}}
\DeclareMathOperator{\Jac}{Jac}

\begin {document}

\newtheorem{thm}{Theorem}
\newtheorem{lem}{Lemma}[section]
\newtheorem{prop}[lem]{Proposition}

\theoremstyle{definition}

\theoremstyle{remark}
\newtheorem{rem}[thm]{Remark}

\title[]{Shifted powers in binary recurrence sequences}
%\author{Samir Siksek}
%\address{Mathematics Institute\\
%	University of Warwick\\
%	Coventry\\
%	CV4 7AL \\
%	United Kingdom}

%\email{s.siksek@warwick.ac.uk}
\author[Michael Bennett]{Michael A. Bennett}
\address{Department of Mathematics, University of British Columbia, Vancouver, B.C., V6T 1Z2 Canada}
\email{bennett@math.ubc.ca}

\author[Sander Dahmen]{Sander R. Dahmen}
\address{Department of Mathematics, VU University Amsterdam, De Boelelaan 1081a, 1081 HV Amsterdam, The Netherlands}
\email{s.r.dahmen@vu.nl}

\author{Maurice Mignotte}
\address{Department of Mathematics, University of Strasbourg, Strasbourg, France}
\email{mignotte@math.u-strasbg.fr}

\author{Samir Siksek}
\address{Mathematics Institute, University of Warwick, Coventry CV4 7AL, United Kingdom}
\email{S.Siksek@warwick.ac.uk}
\thanks{
%Research supported in part by  grants from NSERC, EPSRC, and NWO
The first-named and second-named authors are respectively supported
by NSERC and NWO. The fourth-named author
is supported by an EPSRC Leadership Fellowship EP/G007268/1,
and EPSRC {\em LMF: L-Functions and Modular Forms} Programme Grant
EP/K034383/1.
}

\date{\today}


\begin {abstract}
Let $u_k$ be a Lucas sequence.
A standard technique for determining the perfect powers in the
sequence 
$u_k$ combines bounds coming from linear forms in logarithms
with local information obtained via Frey curves and modularity.
The key to this approach is the fact that the equation $u_k=x^n$
can be translated into a ternary equation of the form $a y^2=b x^{2n}+c$
(with $a$, $b$, $c \in \Z$)
for which Frey curves are available. In this paper we consider
shifted powers in Lucas sequences, and consequently equations of 
the form $u_k=x^n+c$ which do not typically correspond to ternary
equations with rational unknowns. 
However, they do, under certain hypotheses, lead to ternary equations with unknowns in totally real fields,
allowing us to employ Frey curves 
over those fields instead of Frey curves defined over $\Q$. 
We illustrate this approach by showing that
the quaternary Diophantine equation
$x^{2n} \pm 6 x^n+1=8 y^2$ has no solutions in positive integers
$x$, $y$, $n$ with $x$, $n>1$.
\end {abstract}
\maketitle

%------------------------------------------
\section{Introduction} \label{intro}
%------------------------------------------

In \cite{St},  Stewart proved an effective finiteness result for shifted perfect powers in binary recurrence sequences. That is,  if $\{ u_k \}$ is a binary recurrence sequence
for which the equation  
\begin{equation} \label{cur}
x^n+c = u_k
\end{equation}
has a solution in integers $x, n, c$ and $k$, with $n \geq 2$ and $|x| > 1$, then, under mild conditions, $\max \{ |x|, n \}$ is bounded above effectively in terms of $c$ and the recurrence. This statement is actually a consequence of the following more general theorem of Shorey and Stewart \cite{ShSt}.

\begin{thm} (Shorey and Stewart) \label{SS}
Let $a, b, c, d, e$ and $f$ be integers with
$$
(b^2-4ac) (4 acf+bde-ae^2-cd^2-fb^2) \neq 0.
$$
If $x, y$ and $n$ are integers with $|x|>1$ and $n > 2$, satisfying
\begin{equation} \label{dio}
a x^{2n}+b x^n y + c y^2 + d x^n + e y +f = 0,
\end{equation}
then the maximum of $|x|, |y|$ and $n$ is less than a number which is effectively computable in terms of $a, b, c, d, e$ and $f$. Further, if $e^2 \neq 4 cf$ and $x$ and $y$ are integers satisfying
 $$
a x^{4}+b x^2 y + c y^2 + d x^2 + e y +f = 0,
$$
then the maximum of $|x|$ and $|y|$  is less than a number which is effectively computable in terms of $a, b, c, d, e$ and $f$.
\end{thm}

To translate such effective statements to explicit ones regarding equations of the shape (\ref{cur}) or (\ref{dio}) has proven, with current technology, to be a rather challenging problem (and has been accomplished in only a handful of cases -- notably in the determination of perfect powers in the Fibonacci sequence \cite{BMS}). In this paper, we will develop a method which allows us to explicitly find all shifted perfect powers in a  number of classes of Lucas recurrence sequences which are apparently inaccessible to existing techniques in the literature. Our approach combines  lower bounds for linear forms in logarithms (which underlie the proof of Theorem \ref{SS}) with new ideas utilizing connections between Hilbert modular forms and elliptic curves defined over totally real fields.

Whilst we will develop techniques that allow one to carry out such a program in some generality, to focus our exposition we will essentially concentrate on  a single example of an equation of the shape (\ref{dio}), proving  the following.

\begin{thm} \label{main}
The Diophantine equation
\begin{equation} \label{dog}
x^{2n} \pm 6 x^n + 1 = 8y^2
\end{equation}
has no solutions in positive integers $x, n$ and $y$, with $x, n>1$.
\end{thm}

This result is the final ingredient required in work of the first author \cite{Be} on integral points on congruent number curves.
For equation (\ref{dog}), 
it is a fairly routine matter to obtain an absolute bound on $n$  (via Theorem \ref{SS} or otherwise), thereby reducing the problem to that of finding the integral points on a finite collection of hyperelliptic curves. What is much less routine is the approach we take to reduce this bound. Indeed, whilst the problem of determining Fibonacci perfect powers reduces immediately to that of solving ternary equations of the shape
\begin{equation} \label{fibby}
x^2 - 5 y^{2n} = \pm 4,
\end{equation}
for which Frey (or, if you will, Frey--Hellegouarch) curves are immediately available, the fundamental difficulty one encounters in attempting to solve equation (\ref{dog}) is that it is {\it a priori} quaternary rather than ternary. The principal novelty of the paper at hand is that we are able to replace (\ref{dog}) with an equivalent ternary equation over a real quadratic field for which we are able to construct Frey curves which we can, in turn, associate with certain Hilbert modular forms. As in the work of  Bugeaud, Mignotte and Siksek \cite{BMS}, we obtain local information from these Frey curves to  reduce our problem from one of linear forms in three logarithms, to (the computationally more efficient) two logarithms, and subsequently, to find exceptionally strong lower bounds upon $|x|$ for nontrivial solutions to (\ref{dog}), eventually contradicting more general lower bounds for linear forms in (many) complex logarithms.

%-----------------------------------------------------------------------------------------------------
\section{Recurrence sequences : descent to a ternary equation} \label{sec2}
%-----------------------------------------------------------------------------------------------------

To begin the proof of Theorem \ref{main}, let us observe that, if $n=2$, equation (\ref{dog}) with a $(-)$ sign is insoluble modulo $8$; 
with a $(+)$ sign \eqref{dog} defines a genus $1$
curve that is birational over $\Q$ to the rank $0$ elliptic curve with Cremona reference
{\tt 32a1}, and one verifies that the only solutions
on our affine model
satisfy $(|x|,|y|)=(1,1)$.
We may thus suppose that $n>2$ is odd and hence consider the equation
\begin{equation}\label{eqn:main}
x^{2p}+6 x^p+1=8 y^2,
\end{equation}
where $p$ is an odd prime, and $x$ and $y$ are integers.
Note that if we have
$$
u^2+6 u +1=8 y^2,
$$
where $\eta = \mbox{sgn} (u) \in \{ -1, 1 \}$, then $\eta u$ satisfies the recurrence
\begin{equation} \label{rec}
u_{n+1}=6u_n-u_{n-1}+ 12 \, \eta,
\end{equation}
with, say, $(u_0,u_1)=(4-3 \eta,20-3 \eta)$.

Let $K=\Q(\sqrt{2})$ and write $\epsilon=1+\sqrt{2}$ for a fundamental unit of norm $-1$ in $K$. Our main observation that permits application of the so-called modular method is the following.

\begin{lem}\label{lem:descent}
If $(x,y,p)$ is a solution to \eqref{eqn:main} then there exist
integers $k, \ell$ and $s$, and an $\alpha \in \Z[\sqrt{2}]$ such that
\begin{equation}\label{eqn:ppp}
s \, \epsilon^{k} \sqrt{2}-\epsilon^{\ell} \alpha^p=1,
\end{equation}
where $k$ is odd, $s \in \{ -1, 1 \}$,
\begin{equation} \label{stuffy}
\Norm(\alpha)=(-1)^{\ell+1} x \; \; \text{ and } \; \; - \frac{p-1}{2} \leq \ell \leq \frac{p-1}{2}.
\end{equation}
%\margnote{SS: moved restrictions on $l$ into the statement
%of lemma}
\end{lem}
\begin{proof}
We can rewrite \eqref{eqn:main} as 
\[
(x^p+3)^2-8=8 y^2,
\]
whereby $4 \mid x^p+3$ and
\[
y^2 -2 \left(\frac{x^p+3}{4} \right)^2=-1.
\]
Hence
\begin{equation}\label{eqn:rel1}
y+\left(\frac{x^p+3}{4} \right) \sqrt{2} =s \epsilon^{k},
\end{equation}
where $k$ is odd. 

On the other hand, we can also transform equation \eqref{eqn:main} into
\[
\left( \frac{x^p+1}{2} \right)^2 -2 y^2= -x^p
\]
and hence have
\begin{equation}\label{eqn:rel2}
\left( \frac{x^p+1}{2} \right) + y\sqrt{2}= \epsilon^{\ell} \alpha^p,
\end{equation}
where $\alpha$ and $\ell$ satisfy \eqref{stuffy}. 
From \eqref{eqn:rel1}
and \eqref{eqn:rel2}, we deduce \eqref{eqn:ppp}.
\end{proof}

%----------------------------------------------------------------
\section{Linear forms in logarithms} \label{sec3}
%----------------------------------------------------------------

%\margnote{SS: moved the proposition to the beginning
%of the section, to save the lazy reader the trouble
%of reading it!}
The purpose of this section is to prove the following proposition, 
via an appeal to the theory of linear forms in logarithms.
\begin{prop} \label{temp}
If the Diophantine equation
$$
x^{2n} \pm 6 x^n + 1 = 8y^2
$$
has a solution in positive integers $x, n$ and $y$, with $x, n>1$, then $n$ is divisible by an odd prime $p < 2 \cdot 10^{10}$ .
\end{prop}

Either equation (\ref{eqn:ppp}) or (\ref{eqn:rel1}) is a suitable starting point for deriving a linear form in logarithms leading to an absolute upper bound upon $p$; we will appeal to the latter.
Specifically, from (\ref{eqn:rel1}), we have
$$
\frac{\left| x^p+3 \right|}{4} = \frac{\epsilon^{k}+\epsilon^{-k}}{2\sqrt 2}.
$$
%Suppose for now that $k$ is positive, which we may do without loss
%\margnote{MB: replaced $k$ by $|k|$}
%of generality by simply changing the sign of $y$.
It follows that
$$
{|x|^p\over \sqrt 2 \, \epsilon^{|k|}}-1
$$
is ``small'', whereby the same is true of the linear form
\begin{equation} \label{lin-form}
\Lambda = p \log |x| - \log \sqrt 2 - |k| \log \epsilon.
\end{equation}
More precisely, it is easy to verify that
\begin{equation} \label{upper1}
\log |\Lambda | < - p \log |x| +2 .
\end{equation}
%If we put
%$$
%k=q_1b+r, \quad {\rm with}\ |r|<b/2,
%$$
%then we may rewrite $\Lambda$ as
%$$
%\Lambda = b \log\left(p \epsilon^{-q_1} \right)-r\log \epsilon-\log \sqrt 2.
%$$
For any algebraic number $\alpha$ of degree $d$ over $\mathbb{Q}$, we define as usual the {\it absolute logarithmic height} of $\alpha$ by the formula
$$ 
h(\alpha)= \dfrac{1}{d} \left( \log \vert a_{0} \vert + \sum\limits_{i=1}^{d} \log \max \left( 1, \vert \alpha^{(i)}\vert \right) \right), 
$$
where $a_{0}$ is the leading coefficient of the minimal polynomial of $\alpha$ over $\mathbb{Z}$ and the $\alpha^{(i)}$ are the conjugates of $\alpha$ in the field of complex numbers.
The following is the main result (Theorem 2.1) of Matveev \cite{Mat}.
\begin{thm} \label{Matveev} (Matveev)
Let $\mathbb{K}$ be an algebraic number field of degree $D$ over $\mathbb{Q}$ and put $\chi=1$ if $\mathbb{K}$ is real, $\chi=2$ otherwise. Suppose that $\alpha_1, \alpha_2, \ldots, \alpha_n \in \mathbb{K}^*$ with absolute logarithmic heights $h(\alpha_i)$ for $1 \leq i \leq n$, and suppose that
$$
A_i \geq \max \{ D \, h (\alpha_i), \left| \log \alpha_i \right| \}, \; 1 \leq i \leq n,
$$
for some fixed choice of the logarithm. Define
$$
\Lambda = b_1 \log \alpha_1 + \cdots + b_n \log \alpha_n,
$$
where the $b_i$ are integers and  set
$$
B = \max \{ 1, \max \{ |b_i| A_i/A_n \; : \; 1 \leq i \leq n \} \}.
$$
Define, with $e := \exp(1)$, further, 
$$
\Omega =A_1 \cdots A_n, 
$$
$$
C(n) = C(n,\chi) = \frac{16}{n! \chi} e^n (2n+1+2 \chi) (n+2)(4n+4)^{n+1} \left( en/2 \right)^{\chi},
$$
$$
C_0 = \log \left( e^{4.4 n+7} n^{5.5} D^2 \log ( e D) \right) \; \mbox{ and } \; W_0 = \log \left(
1.5 e B D \log (eD) \right).
$$
Then, if $\log \alpha_1, \ldots, \log \alpha_n$ are linearly independent over $\mathbb{Z}$ and $b_n \neq 0$, we have
$$
\log \left| \Lambda \right| > - C(n) \, C_0 \, W_0 \, D^2 \, \Omega.
$$
\end{thm}

We apply this theorem to our situation, with
$$
D=2, \, \chi = 1, \, n=3, \, b_3=p, \, \alpha_3=|x|,
$$
under the assumption that $|x|>1$.
We conclude, after a little work, that
$$
\log |\Lambda| > -  \left( 88626836156 \log p + 232663287513 \right) \log |x|.
$$
Combining this with (\ref{upper1}) (and using that $|x| \geq 7$, an almost immediate consequence of (\ref{eqn:main}) and the supposition that $|x|>1$; succinctly, $2$ is a quadratic residue modulo $x$, whence $x \equiv \pm 1 \pmod{8}$), we 
obtain
the upper bound
\begin{equation} \label{manf}
p < 2{.}772 \cdot 10^{12} =: P_0.
\end{equation}
This bound is the starting point of our analysis. Arguing \`a la Baker, we may thus find an effective absolute upper bound upon $|x|$ in (the finite collection of hyperelliptic) equations  (\ref{eqn:main}). Let us assume, for the remainder of this section, that
\begin{equation} \label{allo}
2 \cdot10^{10} < p < P_0.
\end{equation}
We will show that (\ref{eqn:main}) has no solutions for $p$ satisfying (\ref{allo}), via appeal  to a complicated but slightly sharper lower bound for linear forms in three complex logarithms, due to the third author (Proposition 5.1 of \cite{Mig2}).

\begin{thm} \label{miggy} (Mignotte)
Consider three non-zero  algebraic numbers $\alpha_1$, $\alpha_2$
and $\alpha_3$, which are either all real and ${}>1,$ or all complex of modulus
one and all ${}\not=1$. Further, assume that the three numbers $\alpha_1, \alpha_2$ and $\alpha_3$ are either all multiplicatively independent, or that two of the number are multiplicatively independent and the third is a root of unity.
We also consider three positive
rational integers $b_1$, $b_2$, $b_3$ with $\gcd(b_1,b_2,b_3)=1$, and the linear form
$$
   \Lambda = b_2\log \alpha_2-b_1 \log \alpha_1-b_3\log \alpha_3 ,
$$
where the logarithms of the $\alpha_i$ are arbitrary determinations of the logarithm,
but which are all real or all purely imaginary.
Suppose further that
$$
   b_2 |\log \alpha_2| =
 b_1\,\vert \log \alpha_1 \vert+  b_3 \,\vert\log \alpha_3\vert \pm  \vert\Lambda\vert 
$$
and put
$$
d_1 = \gcd(b_1,b_2) \; \mbox{ and } \;   d_3 = \gcd(b_3,b_2).
$$
Let $\rho\ge e := \exp(1)$  be a real number and set $\lambda = \log \rho$. Let $a_1, a_2$ and $a_3$ be real numbers such that
$$
   a_i \ge  \rho \vert \log \alpha_i \vert
   - \log  \vert \alpha_i\vert +2 D \,{\rm h}\kern .5pt(\alpha_i), \qquad
   i \in \{1, 2, 3 \},
$$
where
$\,D=[\mathbb{Q}(\alpha_1,\alpha_2,\alpha_3) : \mathbb{Q}]\bigm/[\mathbb{R}(\alpha_1,\alpha_2,\alpha_3) : \mathbb{R}]$, and assume further that
$$
\Omega := a_1 a_2 a_3 \geq 2.5 \; \mbox{ and } \; a := \min \{ a_1, a_2, a_3 \} \geq 0.62.
$$
Let $m$ and $L$ be positive integers with $m \geq 3$, $L \geq D+4$  and set
$K = [ m \Omega L ].$
Let $\chi$ be fixed with $0 < \chi \leq 2$ and define
$$
c_1 = \max \{ (\chi m L)^{2/3}, \sqrt{2mL/a} \}, \; 
c_2 = \max \{ 2^{1/3} \, (m L)^{2/3}, L \sqrt{m/a} \}, \;
c_3 = (6 m^2)^{1/3} L,
$$
$$
R_i = \left[ c_i a_2 a_3 \right], \; 
S_i = \left[ c_i a_1 a_3 \right] \; \mbox{ and } T_i = \left[ c_i a_1 a_2 \right],
$$
for $i \in \{ 1, 2, 3 \}$, 
and set
$$
R=R_1+R_2+R_3+1, \; S = S_1+S_2+S_3+1 \; \mbox{ and } \; T = T_1+T_2+T_3+1.
$$
Define
$$
c = \max \left\{ \frac{R}{L a_2 a_3}, \frac{S}{L a_1 a_3}, \frac{T}{L a_1 a_2} \right\}.
$$
Finally, assume that the quantity
$$
\begin{array}{c}
\left( \frac{KL}{2} + \frac{L}{4} - 1 - \frac{2K}{3L} \right) \lambda - (D+1) \log L - 3 g L^2 c \, \Omega  \\ \\
  - D (K-1) \log B - 2 \log K + 2 D \log 1.36
 \end{array}
$$
is positive, 
where
$$
   g={1 \over 4}-{K^2L \over 12RST} \; \mbox{ and } \; 
   B = \frac{e^3 c^2 \Omega^2 L^2}{4K^2 d_1 d_3} \left( \frac{b_1}{a_2}+ \frac{b_2}{a_1} \right) 
    \left( \frac{b_3}{a_2}+ \frac{b_2}{a_3} \right).
$$

\noindent {\bf Then either}
\begin{equation} \label{rups}
\log  \Lambda >  - ( KL + \log ( 3 KL)) \lambda,
\end{equation}
or the following condition holds :

\smallskip
\noindent {\bf either} there exist non-zero rational integers $r_0$ and $s_0$ such that
\begin{equation} \label{rups2}
   r_0b_2=s_0b_1
\end{equation}
with
\begin{equation} \label{rups3}
   |r_0|
   \le \frac{(R_1+1)(T_1+1)}{M-T_1}
    \;  \mbox{ and } \; 
   |s_0| 
   \le \frac{(S_1+1)(T_1+1)}{M-T_1},
\end{equation}
where
$$
M =  \max\Bigl\{R_1+S_1+1,\,S_1+T_1+1,\,R_1+T_1+1,\,\chi \; \tau_1^{1/2} \Bigr\}, \; \; 
\tau_1 = (R_1+1)(S_1+1)(T_1+1),
$$
{\bf or}
there exist rational integers  $r_1$, $s_1$, $t_1$ and $t_2$, with
$r_1s_1\not=0$, such that
\begin{equation} \label{rups4}
   (t_1b_1+r_1b_3)s_1=r_1b_2t_2, \qquad \gcd(r_1, t_1)=\gcd(s_1,t_2 )=1,
\end{equation}
which also satisfy
$$
    |r_1s_1|
   \le \gcd(r_1,s_1) \cdot
   \frac{(R_1+1)(S_1+1)}{M-\max \{ R_1, S_1 \}},
$$
$$
   |s_1t_1| \le \gcd(r_1,s_1) \cdot
   \frac{(S_1+1)(T_1+1)}{M-\max \{ S_1, T_1 \}} 
  $$
and
$$
   |r_1t_2| % &
   \le \gcd(r_1,s_1) \cdot
  \frac{(R_1+1)(T_1+1)}{M- \max \{ R_1, T_1 \}}.
$$
Moreover, when $t_1=0$ we can take $r_1=1$, and
when $t_2=0$ we can take $s_1=1$.
\end{thm}

\medskip

We apply this result with 
$$
b_2=p, \; \alpha_2 = |x|, \;b_1=1, \; \alpha_1 = \sqrt{2}, \; b_3= |k| \; \mbox{ and } \; \alpha_3 =1 + \sqrt{2},
$$
so that we may take
$$
D=2, \; d_1=1, \; d_3 \in \{ 1, p \}, \; a_1 = \frac{\rho+3}{2} \log 2, \; a_2= (\rho+3) \log |x|
$$
and $a_3= (\rho+1) \log ( 1+\sqrt{2})$, whence $a=a_1$. Our goal is to choose $L$, $m,$ $\rho$ and $\chi$ such that (\ref{rups}) contradicts (\ref{upper1}), and (\ref{rups3}) contradicts (\ref{allo}), whereby we necessarily have (\ref{rups4}).

Setting 
$$
L=545, \; m = 25, \, \rho = 5 \; \mbox{ and } \; \chi=2,
$$
we find, after a short Maple computation,  that, for $3 \cdot 10^{10} < p \leq P_0$ and all $|x| \geq 7$, we are in situation (\ref{rups4}),  whereby
there exist integers $r_1, s_1, t_1$ and $t_2$ for which
\begin{equation} \label{pood}
(t_1+r_1 |k|) s_1=r_1t_2 p,
\end{equation}
where
$$
    \left| \frac{r_1s_1}{\gcd(r_1,s_1)} \right|
   \le 80
    \; \mbox{ and } \;
      \left| \frac{s_1t_1}{ \gcd(r_1,s_1)} \right| \le 41.
$$
Similarly, for $2 \cdot 10^{10} < p < 3 \cdot 10^{10}$, we choose 
$$
L=545, \; m = 21, \, \rho = 5 \; \mbox{ and } \; \chi=2,
$$
to deduce the existence of integers $r_1, s_1, t_1$ and $t_2$ with (\ref{pood}) and 
$$
    \left| \frac{r_1s_1}{\gcd(r_1,s_1)} \right|
   \le 75 \; \mbox{ and } \;   \left| \frac{s_1t_1}{ \gcd(r_1,s_1)} \right| \le 39.
$$
Since, in all cases, we assume that $p > 2 \cdot10^{10}$, we thus have
$$
\max \{ |r_1|, |s_1|, |t_1| \} < p.
$$
Hence, from the fact that  $\gcd(r_1, t_1)=\gcd(s_1,t_2 )=1$, it follows that $r_1 = \pm s_1$, whereby
$t_1+r_1 |k| = \pm t_2 p$. Without loss of generality, we may thus write 
$$
u + r |k| = t p,
$$
where $r = |r_1|$ and $t=|t_2|$ are positive integers, $u = \pm t_1$,
$$
|u| \leq 
\left\{
\begin{array}{l}
41 \;  \mbox{ if } \; 3 \cdot 10^{10} < p \leq P_0, \\
39 \;  \mbox{ if }  \; 2 \cdot 10^{10} < p < 3 \cdot 10^{10} \\
\end{array}
\right.
$$
and
$$
r \leq 
\left\{
\begin{array}{l}
80 \;  \mbox{ if } \; 3 \cdot 10^{10} < p \leq P_0, \\
75 \;  \mbox{ if }  \; 2 \cdot 10^{10} < p < 3 \cdot 10^{10} \\
\end{array}
\right..
$$
The linear form $\Lambda$ defined in (\ref{lin-form})  may thus be rewritten as a linear form in two logarithms :
$$
\Lambda =p \log \left( \frac{|x|}{(1+\sqrt{2})^{t/r}} \right) - \log \left( \frac{\sqrt{2}}{(1+\sqrt{2})^{u/r}} \right).
$$
%For future use, it is useful to observe at this juncture that our bounds upon $|u|$ and $r$, together with (\ref{upper1}), imply that we necessarily have
%\begin{equation} \label{hun}
%p \left| \log \left( \frac{|x|}{(1+\sqrt{2})^{t/r}} \right) \right| < 100.
%\end{equation}

We are in position to apply the following state-of-the-art lower bound for linear forms in the logarithms of two algebraic numbers, due to Laurent  (Theorem 2 of \cite{Lau}).  

\begin{thm} \label{laurentlemma} (Laurent)
Let $ \alpha_{1} $ and $ \alpha_{2}$ be multiplicatively independent algebraic numbers, $ h $, $ \rho $ and $ \mu $ be real numbers with $  \rho > 1 $ and $ 1/3 \leq \mu \leq 1$. Set
\begin{center}
$ \begin{array}{ccc}
\sigma=\dfrac{1+2\mu-\mu^{2}}{2}, & \lambda= \sigma \log \rho, & H= \dfrac{h}{\lambda}+ \dfrac{1}{\sigma},
\end{array} $\\
$ \begin{array}{cc}
\omega=2 \left(1+ \sqrt{1+ \dfrac{1}{4H^{2}} } \right), & \theta=\sqrt{1+ \dfrac{1}{4H^{2}} }+ \dfrac{1}{2H}.
\end{array} $
\end{center}
Consider the linear form $ \Lambda=b_{2}\log \alpha_{2}-b_{1}\log \alpha_{1}, $ where $ b_{1} $ and $ b_{2} $ are positive integers.  Put 
$$
 D= \left[  \mathbb{Q}(\alpha_{1}, \alpha_{2} ): \mathbb{Q} \right]/\left[  \mathbb{R}(\alpha_{1}, \alpha_{2} ): \mathbb{R} \right]  
 $$
 and assume that
$$
h \geq \max \left\lbrace D \left( \log \left( \dfrac{b_{1}}{a_{2}}+ \dfrac{b_{2}}{a_{1}} \right) + \log \lambda +1.75 \right)+0.06, \lambda , \dfrac{D \log 2}{2} \right\rbrace ,  
$$

$$
a_{i} \geq \max \left\lbrace 1, \rho \vert \log \alpha_{i} \vert - \log \vert \alpha_{i} \vert + 2Dh(\alpha_{i}) \right\rbrace \ \ \ \ \ (i=1,2), 
$$
and
$$
a_{1}a_{2} \geq \lambda^{2}.
$$

\noindent Then
\begin{equation} \label{laurentall}
 \log \vert \Lambda \vert \geq -C \left( h+ \dfrac{\lambda}{\sigma} \right)^{2} a_{1}a_{2}- \sqrt{\omega \theta} \left(h + \dfrac{\lambda}{\sigma} \right)- \log \left( C' \left(h+ \dfrac{\lambda}{\sigma} \right)^{2} a_{1}a_{2} \right)
 \end{equation}

\noindent with

$$
C=\dfrac{\mu}{\lambda^{3}\sigma} \left( \dfrac{\omega}{6}+ \dfrac{1}{2} \sqrt{\dfrac{\omega^{2}}{9}+ \dfrac{8\lambda \omega^{5/4} \theta^{1/4}}{3 \sqrt{a_{1}a_{2}}H^{1/2} } + \dfrac{4}{3} \left( \dfrac{1}{a_{1}}+ \dfrac{1}{a_{2}} \right) \dfrac{\lambda \omega }{H} } \right)^{2}
$$
and
$$
C'=\sqrt{ \dfrac{C \sigma \omega \theta}{\lambda^{3} \mu} }.
$$
\end{thm}

We apply this result with
$$
b_1=1, \; b_2=p, \; \alpha_1 = \frac{\sqrt{2}}{(1+\sqrt{2})^{u/r}} \; \mbox{ and } \; \alpha_2= \frac{|x|}{(1+\sqrt{2})^{t/r}},
$$
so that $D=2r$,
$$
h (\alpha_1) \leq \frac{\log 2}{2} + \frac{|u|}{2r} \log ( 1+\sqrt{2}) \; \mbox{ and } \;
h (\alpha_2) \leq \log |x| + \frac{t}{2r} \log ( 1+\sqrt{2}).
$$
Further, we may choose
$$
a_1 = \left( 2r + \frac{\rho-1}{2} \right) \log 2 + |u| \left( 2 + \frac{\rho+1}{r} \right) \log ( 1 + \sqrt{2})
$$
and
$$
a_2= \left( 8 \, r + 1 \right) \log |x|.
$$
That this latter choice is a valid one follows from the fact that $|x|^p > (1+\sqrt{2})^{|k|}$, whereby 
$$
t < \left( \frac{u}{|k|}+r \right) \frac{\log |x|}{\log ( 1+ \sqrt{2})},
$$
and the assumption that $\rho < 10^6$, say.

Choosing $\rho=283$ and $\mu = 0.6$, we find, for $p > 3 \cdot 10^{10}$, $|x| \geq	 7$ and all $-41 \leq u \leq 41$, $1 \leq r \leq 80$, that inequality (\ref{laurentall}) contradicts (\ref{upper1}), whilst the same is true (with identical parameter choices), for primes $p$ with $2 \cdot 10^{10} < p < 3 \cdot 10^{10}$, 
$-38 \leq u \leq 38$, $1 \leq r \leq 75$. The final case is when $(r,u)=(75, \pm 39)$ which reduces immediately to that with $(r_0,u_0)=(25,\pm 13)$ upon dividing by the $\gcd (r,u)$. 

Proposition~\ref{temp} thus follows, as desired.

 %---------------------------------------------------------------------------------------
\section{Frey curves and Hilbert modular forms} \label{sec4}
%----------------------------------------------------------------------------------------

%--------------------------------
\subsection{The Frey Curve}
%--------------------------------

We next return to solutions to (\ref{eqn:main}), to which we shall now associate a Frey curve :
\begin{equation}\label{eqn:E}
E_{s,k} \; : \; Y^2=X(X+1)(X+ s\cdot \epsilon^{k} \sqrt{2})
\end{equation}
where the choice of sign $s=\pm 1$ and the value of 
$k$ arises from Lemma~\ref{lem:descent}. 
%It is convenient
%to drop the assumption made at the beginning of the previous
%section that $k$ is positive. 
%\margnote{SS: droping the assumption that $k$ is positive}
By an easy application of Tate's algorithm we find the following.
\begin{lem}
The curve $E_{s,k}$ has minimal discriminant
\[
\Delta_{\mathrm{min}}= 32 \epsilon^{2(k+\ell)} \alpha^{2p}
\]
and conductor
\[
\mathfrak{N}=(\sqrt{2})^9 \cdot \prod_{\fq \mid \alpha} \fq.
\]
\end{lem}

Our goal is to use the arithmetic of this Frey curve to show that any solution to (\ref{eqn:main}) necessarily closely resembles one of  the ``trivial'' ones with $x=1$ and $y=k=\pm 1$.
%\margnote{MB: restated what we mean by trivial solutions}

%--------------------------------
\subsection{Irreducibility}
%--------------------------------

We shall make use of the following result
due to Freitas and Siksek \cite{FS2}, 
which is based on the work of David \cite{DavidI}
and Momose \cite{Momose}.
\begin{prop}\label{prop:irred}
Let $K$ be a totally real Galois number field of degree $d$,
with ring of integers $\OO_K$ and Galois group $G=\Gal(K/\Q)$.
Let $S=\{0,12\}^G$, which we think of as the set of sequences of values $0$, $12$
indexed by $\tau \in G$. %Fix an ordering $\tau_1,\tau_2,\ldots,\tau_d$
For $\bs=(s_\tau) \in S$ and $\alpha \in K$, define the \textbf{twisted norm associated
to $\bs$} by
\[
\cN_\bs(\alpha)= \prod_{\tau \in G} \tau(\alpha)^{s_\tau}.
\]
Let $\epsilon_1,\dots,\epsilon_{d-1}$
be a basis for the unit group of $K$,
and define
\begin{equation}\label{eqn:As}
A_\bs:=\Norm  \left( \gcd ( ( \cN_\bs(\epsilon_1)-1) \OO_K,\ldots, (\cN_\bs(\epsilon_{d-1})-1  ) \OO_K) \right).
\end{equation}
Let $B$ be the least common multiple of the $A_\bs$ taken over all $\bs \ne (0)_{\tau \in G}$,
$(12)_{\tau \in G}$.
Let $p \nmid B$ be a rational prime, unramified in $K$, such that $p \geq 17$ or $p = 11$.
Let $E/K$ be an elliptic curve, and $\fq \nmid p$ be a prime of good reduction for $E$.
Define
\[
P_\fq(X)=X^2-a_\fq(E) X + \Norm(\fq)
\]
to be the characteristic polynomial
of Frobenius for $E$ at $\fq$. Let $r \ge 1$ be an integer such that
$\fq^r$ is principal.
If $E$ is semistable at all $\fp \mid p$
and $\overline{\rho}_{E,p}$ is reducible then
\begin{equation}\label{eqn:res}
p \mid \Res(\, P_\fq(X) \, , \, X^{12 r}-1\, )
\end{equation}
where $\Res$ denotes the resultant of the two polynomials.
\end{prop}

We now return to the case where $E=E_{s,k}$ is our Frey curve \eqref{eqn:E}.
\begin{lem}\label{lem:abirred}
For $E=E_{s,k}$ as above,
$\overline{\rho}_{E,p}$ is irreducible for $p \ge 5$.
\end{lem}
\begin{proof}
We suppose first that $p\ge 17$ or $p=11$, and
apply Proposition~\ref{prop:irred}. The constant $B$
in the proposition is simply
$$
B=\Norm(\epsilon^{12}-1)=-2^5 \cdot 5^2 \cdot 7^2,
$$
whereby if $p \ge 11$ then $p \nmid B$.
Suppose that $\overline{\rho}_{E,p}$ is reducible.
From \eqref{eqn:main}, if $q \equiv 3$, $5 \pmod{8}$
then $q \nmid x$ and so $E$ has good reduction 
at (the inert prime) $q$. 
%with $\mathfrak{a}=(\sqrt{2})^4$
%and hence $A=16$. 

We write $\mu_q$ for the multiplicative order
of $\epsilon$ modulo $q$.
Note that the trace $a_q(E_{s,k})$ depends only on the choice
of sign $s=\pm 1$  and on the value of $k$
modulo $\mu_q$. 
We shall restrict our attention to the following set of primes
\[
Q=\{
3, 5, 11, 13, 19, 29, 43, 59, 83, 109, 131, 139, 251, 269, 307, 419, 461, 659
\}.
\]
The primes $q \in Q$ satisfy $q \equiv 3$, $5 \pmod{8}$ and also
\[
\mu_q \mid 9240=2^3 \cdot 3 \cdot 5 \cdot 7 \cdot 11.
\]
Recall also that $k$ is odd,
and  that $\pm \epsilon^{k}\sqrt{2}-1=\epsilon^{\ell} \alpha^p
\not\equiv 0 \pmod{q}$.
Let
\begin{equation}\label{eqn:S}
S=\{ (t,m) \; : \; \text{$0 \le m < 9240$ odd,\; $t=\pm 1$,\;
$q \nmid (t \cdot \epsilon^m \sqrt{2}-1)$ for all $q \in Q$}
\}.
\end{equation}
It is clear that there is some $(t,m) \in S$ such that
$a_q(E_{s,k})=a_q(E_{t,m})$ for all $q \in Q$.
By Proposition~\ref{prop:irred}, we see that $p$ divides $R_{(t,m)}$ where
\[
R_{(t,m)} = \gcd\{\Res(x^{12}-1, x^2-a_q(E_{t,m})x+q^2) \; :\; q \in Q \}.
\]
Using a short {\tt Magma} \cite{magma} script , we computed $R_{(t,m)}$
for $(t,m) \in S$ and checked that $R_{(t,m)}$ is divisible only by
powers of $2$, $3$, $5$, $7$ and $13$.
Thus $\overline{\rho}_{E,p}$ is
irreducible for $p \ge 17$ or $p=11$.

We now briefly treat $p \in \{ 5, 7, 13 \}$. In each case, we will
in fact show that there is no elliptic curve $E/K$
with full $2$-torsion and a $p$-isogeny. 
Here we found {\tt Magma}'s built-in {\tt Small 
Modular Curve} package invaluable.
\begin{enumerate}
\item[(a)] An elliptic curve $E/K$ with full $2$-torsion and 
a $5$-isogeny is isogenous over $K$ to an elliptic curve
with a $20$-isogeny, and so gives rise to a non-cuspidal
$K$-point on $X_0(20)$. A model for $X_0(20)$ is given by
the elliptic curve
\[
X_0(20) \; : \; y^2 = x^3 + x^2 + 4x + 4,
\]
which has Cremona reference 20A1.
This has rank $0$ over $K$, and in fact a full list
of $K$-points is $\{\infty, (4 , \pm 10), (0 , \pm 2), (-1 , 0 )\}$.
These points are all cuspidal, completing the proof for $p=5$.
\item[(b)] An elliptic curve $E/K$ with a $7$-isogeny
and a point of order $2$ gives rise to a non-cuspidal
$K$-point on $X_0(14)$. A model for $X_0(14)$ is
given by the elliptic curve
\[
X_0(14) \; : \; y^2 + x y + y = x^3 + 4 x - 6,
\] 
which has Cremona reference 14A1.
This again has rank $0$ over $K$. The $K$-points
are 
$\{\infty, (9 , 23), 
(1 , -1 ), (2 , -5 ), 
(9 , -33 ),
(2 , 2 ) 
\}$. The first four points are cusps. The last two
correspond to elliptic curves with $j$-invariants
$-3375$ and $16581375$. It turns out that elliptic curves
with these $j$-invariants have only one point of order $2$ over $K$.
This completes the proof for $p=7$.
\item[(c)] An elliptic curve $E/K$ with a $13$-isogeny
and a point of order $2$ gives rise to a non-cuspidal
$K$-point on $X_0(26)$, which has genus $2$. We shall in fact work
with $X_0(26)/\langle w_{13} \rangle$, where $w_{13}$ is the 
Atkin-Lehner involution. This quotient is the
elliptic curve with Cremona reference 26B1:
\[
X_0(26)/\langle w_{13} \rangle \; : \; y^2 + x y + y = x^3 - x^2 - 3 x + 3.
\] 
Again this has rank $0$ over $K$, and a full list of $K$-points
is given by $\{ 
\infty, 
(-1 , -2), 
(-1 , 2 ),
(1 , -2 ), 
(1 , 0 ), 
(3 , -6 ), 
(3 , 2 )
\}$. The only $K$-points we obtain by pulling back to $X_0(26)$ 
are cusps. This completes the proof.
\end{enumerate}
\end{proof}

%--------------------------------
\subsection{Level-lowering}
%--------------------------------

Suppose that $p \ge 5$.
%From the fact that $3$ is inert in $K$ and $E=E_{s,k}$ has good reduction at $3\cdot \Z[\sqrt{2}]$,
%we know that $E$ is modular \cite[Theorem 4.3]{FreitasRecipes}. 
By \cite{FLS}, elliptic curves 
over real quadratic
fields are modular; in particular $E$ is modular.
We now
apply the standard level-lowering recipe found in \cite{FS}
and based on theorems of Fujiwara, Jarvis and Rajaei (it is
here that we require Lemma~\ref{lem:abirred}). From the recipe
we know that $\overline{\rho}_{E,p} \sim \overline{\rho}_{f,\fp}$
for some Hilbert newform over $K$ of level $\mathfrak{M}=(\sqrt{2})^9$
and prime ideal $\fp \mid p$. Using {\tt Magma}, 
we find that the space of Hilbert newforms of level $\mathfrak{M}$
is $8$-dimensional, and in fact decomposes into $8$ rational
eigenforms. 

Through a small search, we found $8$ elliptic curves over $K$ of conductor 
$\mathfrak{M}$. By computing their traces at small prime ideals,
we checked that they are in fact pairwise non-isogenous. These
elliptic curves are all modular by the same theorem cited above, and
hence must correspond to the $8$ Hilbert newforms of level
$\mathfrak{M}$. Thus $\overline{\rho}_{E,p} \sim \overline{\rho}_{F_i,p}$
where $F_1,\dots,F_8$ are the $8$ elliptic curves given as follows :
\begin{gather*}
F_1 \; : \; 
Y^2 = X^3 + \sqrt{2} X^2 + (\sqrt{2} - 1) X,\\
F_2 \; : \;  
Y^2 = X^3 + (-\sqrt{2} + 3) X^2 + (-\sqrt{2} + 2) X,\\
F_3 \; : \; 
 Y^2 = X^3 + (2 \sqrt{2} - 1) X^2 + (-\sqrt{2} + 2) X,\\
F_4 \; : \; 
 Y^2 = X^3 + (\sqrt{2} - 2) X^2 + (-\sqrt{2} + 1) X,\\
F_5 \; : \; 
 Y^2 = X^3 + (-\sqrt{2} + 1) X^2 - \sqrt{2} X,\\
F_6 \; : \; 
 Y^2 = X^3 + (\sqrt{2} - 1) X^2 - \sqrt{2} X,\\
F_7 \; : \; 
 Y^2 = X^3 + (\sqrt{2} + 3) X^2 + (\sqrt{2} + 2) X,\\
F_8 \; : \; 
 Y^2 = X^3 - \sqrt{2} X^2 + (-\sqrt{2} - 1) X.
\end{gather*}

\begin{lem}\label{lem:KO}
Let $E=E_{s,k}$ and let $F$ be one of the eight 
elliptic curves $F_1,\dots,F_8$ above.
Suppose that $\overline{\rho}_{E,p} \sim \overline{\rho}_{F,p}$ and
let $\fq \nmid 2$ be a prime ideal of $K$, and 
$q$ be the rational prime such that $\fq \mid q$. 
\begin{enumerate}
\item[(i)] If $\fq \nmid (s \epsilon^k \sqrt{2}-1)$ then
$a_\fq(E) \equiv a_\fq(F) \pmod{p}$.
\item[(ii)] If $\fq \mid (s \epsilon^k \sqrt{2}-1)$ 
and $q \not\equiv 7 \pmod{8}$ then
$\Norm(\fq)+1 \equiv a_\fq(F) \pmod{p}$.
\item[(iii)] If $\fq \mid (s \epsilon^k \sqrt{2}-1)$ 
and $q \equiv 7 \pmod{8}$ then
$\Norm(\fq)+1 \equiv -a_\fq(F) \pmod{p}$.
\end{enumerate}
\end{lem}
\begin{proof}
From the model \eqref{eqn:E},
it is easy to see that $E_{s,k}$ has good reduction at $\fq$ in case (i),
split multiplicative reduction in case (ii) and non-split
multiplicative reduction in case (iii).
If $q \ne p$,
the lemma follows by comparing traces of the images of Frobenius at $\fq$
for the representations $\overline{\rho}_{E,p}$ and 
$\overline{\rho}_{F,p}$. For $q=p$ see \cite[Proposition 3]{KO}.
\end{proof}

\begin{lem}\label{lem:sign}
The $s=\pm 1$ sign in \eqref{eqn:ppp}, \eqref{eqn:rel1} and \eqref{eqn:E}
 is in fact $+1$. Moreover, either
$k \equiv -1 \pmod{9240}$ and 
$\overline{\rho}_{E,p} \sim \overline{\rho}_{F_2,p}$
or $k \equiv 1 \pmod{9240}$ and
$\overline{\rho}_{E,p} \sim \overline{\rho}_{F_7,p}$.
\end{lem}
\begin{proof}
We shall use ideas from the proof of Lemma~\ref{lem:abirred}.
For 
a prime ideal $\fq$ we let $\mu_\fq$ be the order of $\epsilon$
modulo $\fq$. 
In the proof of Lemma~\ref{lem:abirred}, we restricted ourselves
to inert primes. It is useful to now also use split primes.
Let
\[
\mathfrak{Q}=\{ \fq \; : \; \text{$\mu_\fq \mid 9240$, \;
$3 \le \Norm(\fq) < 1000$} \}.
\]

Fix $i \in \{ 1, 2, \dots, 8 \}$ and suppose $\overline{\rho}_{E,p} \sim
\overline{\rho}_{F_i,p}$ where $E=E_{s,k}$. 
We let $S$ be given by \eqref{eqn:S}, whereby we know that there
is some $(t,m) \in S$ such that $s=t$ and $k \equiv m \pmod{9240}$.
Let $\fq \in \mathfrak{Q}$, and let $q$ be the
rational prime satisfying $\fq \mid q$.
By Lemma~\ref{lem:KO} we have the following. 
\begin{enumerate}
\item[(a)] If $\fq \nmid (t \epsilon^m \sqrt{2}-1)$ then
$a_{\fq}(E_{t,m})=a_{\fq}(E_{s,k}) \equiv a_{\fq}(F_i) \pmod{p}$.
\item[(b)] If $\fq \mid (t \epsilon^m \sqrt{2}-1)$ 
and $q \not\equiv 7 \pmod{8}$ then
$\Norm(\fq)+1 \equiv a_\fq(F_i) \pmod{p}$.
\item[(c)] If $\fq \mid (t \epsilon^m \sqrt{2}-1)$ 
and $q \equiv 7 \pmod{8}$ then
$\Norm(\fq)+1 \equiv -a_\fq(F_i) \pmod{p}$.
\end{enumerate}

Define
\[
\beta_{t,m,i}(\fq)=\begin{cases}
a_{\fq}(E_{t,m}) -a_{\fq}(F_i) & \text{if $\fq \nmid (t \epsilon^m \sqrt{2}-1)$}\\
\Norm(\fq)+1 - a_\fq(F_i) &
\text{
if $\fq \mid (t \epsilon^m \sqrt{2}-1)$ and $q \not\equiv 7 \pmod{8}$
}
\\
\Norm(\fq)+1 + a_\fq(F_i) &
\text{
if $\fq \mid (t \epsilon^m \sqrt{2}-1)$
and $q \equiv 7 \pmod{8}$
}.
\end{cases}
\]
We let 
\[
\gamma_{t,m,i}= \gcd ( \beta_{t,m,i}(\fq) : \fq \in \mathfrak{Q}
).
\]
Note that if $s=t$, $k \equiv m \pmod{9240}$ and $\rho_{E,p} \sim \rho_{F_i,p}$, 
then
$p \mid \gamma_{t,m,i}$. Using a {\tt Magma} script, we computed 
$\gamma_{t,m,i}$ for all $(t,m) \in S$ and $i \in \{ 1,2, \dots, 8 \}$.
We found that all of these are divisible only by $2$ and $3$,
except for $\gamma_{1,9239,2}$ and $\gamma_{1,1,7}$
which are both zero. Hence $s=+1$, and either
$k \equiv -1 \pmod{9240}$ and $F_i=F_2$
or $k \equiv 1 \pmod{9240}$ and $F_i=F_7$.
\end{proof}

Note, in fact, that $F_2$ is isomorphic to $E_{1,-1}$ and
$F_7$ is isomorphic to $E_{1,1}$ (which explains why
$\gamma_{1,9239,2}=\gamma_{1,1,7}=0$).

We now simplify our Frey curve in \eqref{eqn:E} to take
account of the sign $s=+1$. We denote the new Frey curve
by
\[
E_k \; : \; Y^2=X(X+1)(X+\epsilon^k \sqrt{2}).
\]

Equations \eqref{eqn:main} and \eqref{eqn:ppp}
satisfy the following useful symmetry.
%\margnote{SS: changed statement (and proof) of lemma
%to make symmetry of solutions clearer.}
\begin{lem}\label{lem:choice}
Let $(x,y,k,\ell,\alpha)$ be a solution to equations  \eqref{eqn:main} and
\eqref{eqn:ppp}, satisfying \eqref{stuffy}.
Then $(x,-y,-k,-\ell,(-1)^\ell \overline{\alpha})$ is also
a solution to equations \eqref{eqn:main} and \eqref{eqn:ppp}, also satisfying \eqref{stuffy}. 
\end{lem}
\begin{proof}
The lemma follows on conjugating equations \eqref{eqn:ppp} and
\eqref{eqn:rel1},  
observing 
that $\overline{\epsilon}=-\epsilon^{-1}$, and 
recalling that $k$ is odd.
\end{proof}
%\begin{lem}\label{lem:choice}
%Let $F=F_7$.
%After possibly changing the sign of $y$ in \eqref{eqn:main},
%we have $k \equiv 1 \pmod{9240}$ and $\overline{\rho}_{E_k,p}
%\sim \overline{\rho}_{F,p}$.
%\end{lem}
%\begin{proof}
%By replacing $y$ by $-y$ in \eqref{eqn:main}, we replace 
%the term $\epsilon^k$ in \eqref{eqn:rel1} by
%\[
%-\left(\overline{\epsilon}\right)^k=-(-\epsilon^{-1})^k=\epsilon^{-k}
%\]
%since $k$ is odd. In other words, replacing $y$ by $-y$
%leads to replacing $k$ by $-k$.
%\end{proof}

%---------------------------------------------------------------
\section{Arithmetic information from Frey curves}
%---------------------------------------------------------------

We continue with the same notation as in the previous section. Our basic goal is to obtain lower bounds for the exponent $|k|$ in the event that $k \neq 1$.

\subsection{Sieving : part I}

We begin by sharpening (part of) Lemma \ref{lem:sign}.
\begin{prop}\label{prop:1modM}
Suppose that $p \ge 19$ and let
\begin{equation}\label{eqn:M}
M=9240  \prod_{\overset{3 \le \ell < 2.4 \times 10^5}{\text{$\ell$ prime}}} \ell.
\end{equation}
%\margnote{SS: clarified the lemma
%by allowing $\pm 1$ and made
%appropriate change to the proof}
Then $k \equiv \pm 1 \pmod{M}$. In particular, 
$k=\pm 1$ or $\log_{10}{\lvert k \rvert} \ge 103944$.
\end{prop}
\begin{proof}
We start by defining $M_0=9240=2^3 \cdot 3 \cdot 5 \cdot 7 \cdot 11$.
We know from Lemma~\ref{lem:sign} that $k \equiv \pm 1 \pmod{M_0}$.
By Lemma~\ref{lem:choice}, we may assume that $k \equiv 1 \pmod{M_0}$,
and will use this to deduce that $k \equiv 1 \pmod{M}$.
It then plainly follows that if $k \equiv -1 \pmod{M_0}$
then $k \equiv -1 \pmod{M}$. 

Suppose $k \equiv 1 \pmod{M_0}$. Let 
\[
\ell_1=3, \quad \ell_2=5, \quad \ell_3=7, \dots
\]
be the sequence of primes starting with $3$. We define 
$M_{n}=\ell_n \cdot M_{n-1}$.  We will show inductively that
$k \equiv 1 \pmod{M_{n}}$ until $M_n=M$.
A direct computation, showing that $M$ somewhat exceeds $10^{103944}$, yields the last
statement in the proposition. 

For the inductive step, suppose $k \equiv 1 \pmod{M_{n-1}}$.
Our strategy is to write down a small set 
$\mathfrak{Q}$
of odd prime ideals
$\fq$ of $K$ satisfying
\[
\ell_n \mid \mu_\fq, \qquad \mu_\fq \mid M_{n};
\] 
here as before, 
$\mu_\fq$ is the multiplicative order of $\epsilon$ modulo $\fq$.
Let
\[
\mathcal{K}=\{1,1+M_{n-1},1+2M_{n-1},\dots,1+(\ell_{n}-1) M_{n-1} \}.
\]
We know that $k \equiv m \pmod{M_n}$ for some $m \in \mathcal{K}$.
By the previous section, $p$ divides $\beta_{1,m,7}(\fq)$ for all
$\fq \in \mathfrak{Q}$. For $m \in \mathcal{K}$,
we compute $\gcd_{\fq \in \mathfrak{Q}} (\beta_{1,m,7}(\fq))$. With a 
sufficiently large enough initial set $\mathfrak{Q}$, we
found that this gcd is divisible only by primes $\le 17$ unless
$m=1$ (in which case it is $0$). This shows that $k \equiv 1 \pmod{M_{n}}$.
The {\tt Magma} script executing this proof took roughly 218 hours to run on a 2200MHz AMD Opteron.
\end{proof}

\subsection{Sieving : part II}

In view of Proposition~\ref{prop:1modM}, we suppose $k \equiv \pm 1 \pmod{M}$
where $M$ is given by \eqref{eqn:M}. The objective of this subsection is to
show the following.  
%\margnote{SS: changed statement of prop to allow $\pm$ signs,
%and adapted the proof.}
\begin{prop}\label{prop:k1modp}
For primes $19 \leq p < 2 \cdot10^{10}$ we have $k \equiv \ell \equiv \pm 1 \pmod{p}$,
where $k$ and $\ell$ are the exponents in \eqref{eqn:ppp}. In particular, via 
\eqref{stuffy}, we have that $\ell=\pm 1$.
\end{prop}

We shall suppose that $k \equiv 1 \pmod{M}$ and deduce $k \equiv \ell \equiv 1
\pmod{p}$. The Proposition then follows from Lemma~\ref{lem:choice}.
Fix $p \ge 19$. Inspired by \cite[Lemma 7.4]{BMS},
we choose an auxiliary integer $q$ satisfying certain conditions, which we
enumerate as needed.
The first two conditions are the following.
\begin{enumerate}
\item[(i)] $q \equiv 1 \pmod{8}$ is prime.
\item[(ii)] $q =np+ 1$, where $n$ is an integer.
\end{enumerate}
Fix $\delta \in \F_q$ satisfying $\delta^2 \equiv 2 \pmod q$ (which we may do, via assumption  (i)). Let
$$
\fq_1=q \OO_K + (\sqrt{2}-\delta)\OO_K \; \mbox{ and } \; 
\fq_2=q \OO_K + (\sqrt{2}+\delta)\OO_K.
$$
Then $\fq_1$ and $\fq_2$ are prime ideals with residue field $\F_q$, 
\[
\sqrt{2} \equiv \delta \pmod{\fq_1} \; \mbox{ and } \;  \sqrt{2} \equiv -\delta \pmod{\fq_2}.
\]
Let
$$
\mathcal{F}_1 \; : \; Y^2=X^3+(\delta+3) X^2+(\delta+2)X
$$
and 
$$
\mathcal{F}_2 \; : \; Y^2=X^3+(-\delta+3) X^2+(-\delta+2)X.
$$
These are the reductions of $F_7$ modulo $\fq_1$ and $\fq_2$,
respectively. We shall suppose that $q$ is chosen so that
 the following condition on the traces of $\mathcal{F}_i$ is satisfied.
\begin{enumerate}
\item[(iii)] $a_q(\mathcal{F}_i) \not\equiv \pm 2 \pmod{p}$
for $i=1$, $2$.
\end{enumerate}
Note that $q \equiv 1 \pmod{p}$ by condition (ii). By Lemma~\ref{lem:KO}, we
see that 
$$
\fq_i \nmid (\epsilon^k \sqrt{2}-1). 
$$
Thus by 
Lemma~\ref{lem:descent}, $q \nmid x$. Appealing to Lemma~\ref{lem:KO}
again, we have
\begin{equation}\label{eqn:traces}
a_{\fq_1}(E_k) \equiv a_q(\mathcal{F}_1) \pmod{p} \; \mbox{ and } \; 
a_{\fq_2}(E_k) \equiv a_q(\mathcal{F}_2) \pmod{p}.
\end{equation}

By \eqref{eqn:rel1} and Lemma~\ref{lem:sign}, we have
\[
y+\left(\frac{x^p+3}{4} \right) \sqrt{2} =\epsilon^{k},
\]
whence, using the fact that $k$ is odd,
\begin{equation}\label{eqn:xepsk}
\left(\frac{x^p+3}{2} \right) \sqrt{2}= \epsilon^k+\epsilon^{-k}.
\end{equation}
Let
\[
\mu_n(\F_q)=\{ \mu \in \F_q^* \; : \; \mu^n \equiv 1 \pmod{q} \},
\]
whereby, as $q \nmid x$ and $q=np+1$, we see that
\[
(x^p \mod{q}) \in \mu_n(\F_q).
\]
Let
\[
\mathcal{W}_q=\left\{ e \in \F_q \; :
\;  e+e^{-1}=\frac{(\mu+3)\delta}{2} \; \text{for some  $\mu \in \mu_n(\F_q)$}  \right\}.
\]
By \eqref{eqn:xepsk}, we have that 
\begin{equation}\label{eqn:epske}
\epsilon^k \equiv e \pmod{\fq_1} \; \mbox{ and } \;  \epsilon^{k} \equiv -e^{-1} \pmod{\fq_2},
\end{equation}
for some $e \in \mathcal{W}_q$. Moreover, since $\fq_i \nmid (\epsilon^k \sqrt{2}-1)$, if we define
\[
\mathcal{X}_q=\{e \in \mathcal{W}_q \; : \; e \not\equiv \delta^{\pm 1} \pmod{q}\},
\]
then necessarily \eqref{eqn:xepsk} holds for some $e \in \mathcal{X}_q$.
In practice, we hope to be able to find a prime $q$ satisfying the foregoing
and forthcoming conditions, with $n$ small. Computing $\mathcal{W}_q$
amounts to solving $n$ quadratic equations in $\F_q$. The size of
$\mathcal{W}_q$ and thus $\mathcal{X}_q$ is at most $2n$. Thus 
we know that the reduction
of $\epsilon^k$ modulo $\fq_1$ belongs to this relatively small set.
We will now refine $\mathcal{X}_q$, defining subsets 
$\mathcal{Y}_q$, and $\mathcal{Z}_q$ that also contain
the reduction of $\epsilon^k$ modulo $\fq_1$.
To do this,
also suppose that $q$ is chosen so that the following condition is satisfied.
\begin{enumerate}
\item[(iv)] $n \mid M$.
\end{enumerate}
Since $k \equiv 1 \pmod{M}$, we see that $(q-1) \mid (k-1) p$, whereby
\[
(\epsilon^{k-1})^p \equiv 1 \pmod{\fq_i}.
\]
Note that $\epsilon \equiv 1+\delta \pmod{\fq_1}$.
Let
\[
\mathcal{Y}_q=\left\{  e \in \mathcal{X}_q \; : \; \left(\frac{e}{1+\delta}\right)^p \equiv 1 \pmod{q}\right\}.
\]
We see that \eqref{eqn:epske} holds for some $e \in \mathcal{Y}_q$.
Heuristically, the probability that a random element of $\mathcal{X}_q$
belongs to $\mathcal{Y}_q$ is $1/n$. Since $\# \mathcal{X}_q \le 2n$, 
we expect that $\#\mathcal{Y}_q=O(1)$. 

For the next refinement,
we will use information derived from the modular approach,
as given in \eqref{eqn:traces}.
Write
$$
\mathcal{E}_{1,e} \; : \;
Y^2=X(X+1)(X+e \cdot \delta)
\; \mbox{ and } \; 
\mathcal{E}_{2,e} \; : \;
Y^2=X(X+1)(X+e^{-1}\cdot \delta).
$$
Given \eqref{eqn:epske}, 
these two elliptic curves over $\F_q$ 
are the reductions of the Frey curve $E_k$ modulo $\fq_1$, $\fq_2$.
Let $\mathcal{Z}_q$ be the set of $e \in \mathcal{Y}_q$ such that
\[
a_q(\mathcal{E}_{1,e}) \equiv a_q(\mathcal{F}_1) \pmod{p} \; \mbox{ and } \; 
a_q(\mathcal{E}_{2,e}) \equiv a_q(\mathcal{F}_2) \pmod{p}.
\]
We know from \eqref{eqn:traces} that \eqref{eqn:epske} holds 
for some $e \in \mathcal{Z}_q$. Note that $\mathcal{Z}_q$
cannot be empty, as the value $k=1$ leads to a solution to our 
original equation~\eqref{eqn:main}. Thus certainly,
$1+\delta$ (which is the reduction of $\epsilon$ modulo $\fq_1$)
must appear in $\mathcal{Z}_q$. This, of course, is a useful check
on the correctness of our computations.

It is reasonable to expect on probabilistic grounds that
$\mathcal{Z}_q=\{1+\delta\}$. In fact, this is one of our two final
assumptions on $q$. 
\begin{enumerate}
\item[(v)] $\mathcal{Z}_q=\{ 1+\delta\}$.
\item[(vi)] $(1+\delta)^n \not \equiv 1 \pmod{q}$.
\end{enumerate}
From (v), we see that
\[
(1+\delta)^k \equiv \epsilon^k \equiv 1+\delta \pmod{\fq_1}.
\]
Thus the multiplicative order of $1+\delta$ in $\F_q^*$ divides
$k-1$. Since $q-1=np$, we have from (vi) that
this multiplicative order must be divisible by $p$. Therefore,
$k \equiv 1 \pmod{p}$. 
We now turn our attention to $\ell$. 
Reducing \eqref{eqn:ppp} modulo $\fq_1$
(and recalling the value of the sign from Lemma~\ref{lem:sign})
we have
\[
(1+\delta)\delta - (1+\delta)^\ell \alpha^p \equiv 1 \pmod{\fq_1}.
\]
As $\delta^2 \equiv 2 \pmod{\fq_1}$, it follows that 
\[
(1+\delta)^\ell \alpha^p \equiv 1+\delta \pmod{\fq_1}
\]
and so
\[
(1+\delta)^{n(\ell-1)} \equiv \alpha^{-np} \equiv 1 \pmod{\fq_1}.
\]
From (vi), we conclude that $\ell \equiv 1 \pmod{p}$.
The following lemma summarizes the above.
%\margnote{SS: Added assumption that $k \equiv 1 \pmod{M}$
%in the statement of the lemma.}
\begin{lem}
Suppose $k \equiv 1 \pmod{M}$.
Suppose there is a prime $q$ satisfying conditions (i)--(vi).
Then 
$$
k \equiv \ell \equiv 1 \pmod{p}.
$$
\end{lem}

\begin{proof}[Proof of Proposition~\ref{prop:k1modp}]
As observed previously, it is sufficient to suppose
%\margnote{SS: Proof adjusted.}
that $k \equiv 1 \pmod{M}$ and show that 
$k \equiv \ell \equiv 1 \pmod{p}$, for primes $p$
in the range $19 \le p < 2\cdot 10^{10}$.
We used a {\tt Magma} script which for each prime $p$ in this range,
finds a prime $q$ satisfying conditions (i)--(vi) above.
The total processor time for the proof is roughly 1946 hours,
although the computation,
running on a 2200MHz AMD Opteron,
 was spread over 10 processors,
making the actual computation time less than nine days. 
\end{proof}

%---------------------------------------------------------
\section{Reducing our upper bound on $p$}
%----------------------------------------------------------

From Proposition \ref{temp}, we may suppose that $p < 2 \cdot 10^{10}$. The goal of this section is to reduce this upper bound still further. 
Let $p\geq 19$. We have, according to Lemma~\ref{lem:sign}, that $s=1$ in~\eqref{eqn:ppp}.
We will assume for the remainder of this section that $k > 1$.
Appealing to Proposition \ref{prop:k1modp}, there exists a positive integer $k_0$  such that $k=p k_0 \pm 1$, whereby we may now rewrite (\ref{eqn:ppp}) as one of
\begin{equation} \label{eqn:new} 
\alpha^p - \sqrt{2} \left( \epsilon^{k_0} \right)^p = \pm 1 - \sqrt{2},
\end{equation}
so that
\begin{equation} \label{wombat}
0 < \Lambda_1 = \log ( \sqrt{2} ) - p \log \left( \frac{\alpha}{\epsilon^{k_0}} \right) < 
\epsilon \cdot  \alpha^{-p}.
\end{equation}

Applying Theorem \ref{laurentlemma}, with
$$
b_2=1, \; \alpha_2 = \sqrt{2}, \; b_1=p, \; \alpha_1 = \alpha/\epsilon^{k_0}  \; \mbox{ and } \; D=2,
$$
we may take
$$
a_1 = (\rho-1) \frac{\log 2}{ 2p}  + 2 \log ( \alpha) + 4 \, h (\alpha)  \leq  (\rho-1) \frac{\log 2}{ 2p}  + 6 \log ( \alpha)
$$
and
$$ 
a_2 = \left( \frac{\rho+3}{2} \right) \log 2.
$$
We choose
$\rho=27$ and $\mu = 1/3$, whereby a short calculation ensures that inequality (\ref{laurentall}), together with the fact that   $\alpha \geq \sqrt{7}$, contradicts (\ref{wombat}) for $p > 1637$. Computing the first $10^7$ terms in recursion (\ref{rec}), verifying that none of them coincide with solutions to (\ref{eqn:main}) and noting that the $10^7$-th term exceeds $\exp(8.8 \cdot 10^6)$, we thus have $|x| > 10^{2348}$ and hence $|\alpha| > 10^{1174}$.
Now choosing $\rho=31$ and $\mu = 1/3$ in Theorem \ref{laurentlemma}, and using our new lower bound upon $\alpha$, we find that $p \leq 941$.

%If $\rho$ is suitable large, we may take $h=2 \log p$.

%------------------------------------------------------------------------
\subsection{Handling the values of $19 \leq p \leq 941$}
%------------------------------------------------------------------------

From equation (\ref{eqn:new}), it suffices to solve the family of Thue equations
\begin{equation} \label{gut}
X^p -  \sqrt{2} \, Y^p = \pm 1 - \sqrt{2},
\end{equation}
in integers $X, Y \in \mathbb{Z} [ \sqrt{2}]$, for primes $p$ with $19 \leq p \leq 941$, to handle~\eqref{eqn:main} for these values of $p$. We will do this under the additional assumption that
$$
Y = \epsilon^{k_0} > 1.
$$

We argue as in the proof of Proposition 9.3 of \cite{BMS}. Define $\omega =
2^{1/2p}$ to be real and positive and let $\zeta$ be a primitive $p$-th root of
unity. Set $K = \mathbb{Q}(\omega)$ and $L=\mathbb{Q}(\omega,\zeta)$. It
is straightforward to check that an integral basis for $K$
is $1,\omega,\ldots,\omega^{2p-1}$, and to use this
%\margnote{SS: removed the reference to Magma. It
%must be using a theorem instead of doing the computation.
%I don't know what the theorem is, but it's easy to check
%the discriminant with pen and paper!}
to deduce that the discriminant of $K$ is 
$D_K = 2^{4p-1} p^{2p}$. Moreover, the unit rank of $K$ is
$p$ and its Galois closure is $L$. We write
$\epsilon_{1,1}, \ldots, \epsilon_{1,p}$ for a system of fundamental units  of
$K$ which, via Lemma 9.6 of \cite{BMS}, we may suppose to satisfy
\begin{equation} \label{ilk}
\prod_{i=1}^p h ( \epsilon_{1,i}) \leq 2^{1-p} (p!)^2 (2p)^{-p} R_K,
\end{equation}
where $R_K$ denotes the regulator of $K$. Further,  the absolute values of the inverse of the regulator matrix corresponding to $\epsilon_{1,1}, \ldots, \epsilon_{1,p}$ are bounded above by $(p!)^2 2^{-p} \log^3 (6p)$. There thus exist integers $b_1, \ldots, b_p$ such that
$$
X - \omega Y = \pm \epsilon_{1,1}^{b_1} \cdots \epsilon_{1,p}^{b_p},
$$ 
where
$$
B= \max \{ |b_1|, \ldots, |b_p| \} \leq 2^{1-p} p (p!)^2 \log^3 (6p) \, h(X-\omega Y).
$$
We will assume that $k$ is positive; the argument for negative $k$ is similar and leads to an identical conclusion. From (\ref{gut}), considering imaginary parts, we have that
$$
\left| X - \omega Y \right| = ( \sqrt{2} \pm 1) \prod_{i=1}^{p-1} \left| X - \omega \zeta^i Y \right|^{-1}
\leq ( \sqrt{2}+1) \, \omega^{1-p} \, |Y|^{1-p} \, \prod_{i=1}^{(p-1)/2} \sin^{-2} (\pi i/p),
$$
whence
\begin{equation} \label{lush}
\left| X - \omega Y \right|< 2^p \, |Y|^{1-p}.
\end{equation}
It follows that
$$
 \left| X - \omega \zeta^j Y \right|  \leq \omega |Y| \left| 1 - \zeta^j \right| + 2^p \, |Y|^{1-p} < 2.1 \, |Y|,
 $$
 where, from the assumption that $k > 1$, we have
 \begin{equation} \label{hoop}
\log |Y| = \frac{k \pm 1}{p} \, \log (1+\sqrt{2}) \geq \frac{M}{p} \log (1+ \sqrt{2}) > 10^{103942} \log (1+ \sqrt{2}).
\end{equation}
We therefore estimate that
$$
h(X-\omega Y) < 0.75 + \log |Y| < 1.01 \, \log |Y|,
$$
whereby, crudely,
\begin{equation} \label{why}
B < p^{2p} \, \log |Y|.
\end{equation}

Define $\epsilon_{2,i}$ and $\epsilon_{3,i}$ for $1 \leq i \leq p$ as the images under $\sigma$ and $\sigma^2$ of $\epsilon_{1,i}$,  where $\sigma$ is the field automorphism that sends $\omega$ to $\omega \zeta$ and fixes $\zeta$. If, following  Siegel, we set
$$
\lambda = \frac{-\zeta}{1+\zeta} \cdot \frac{X-\omega Y}{X - \zeta \omega Y} =
 \frac{X-\zeta^2 \omega Y}{X - \zeta \omega Y} \cdot \frac{1}{1+\zeta} - 1,
$$
then we have
$$
\lambda = \left( \frac{\epsilon_{3,1}}{\epsilon_{2,1}} \right)^{b_1}  \cdots 
 \left( \frac{\epsilon_{3,p}}{\epsilon_{2,p}} \right)^{b_p}  \frac{1}{1+\zeta} - 1.
$$
Since $|\zeta/(1+\zeta)| < 1$, arguing crudely, from (\ref{lush}), we have that $|\lambda| < 2^{p} \, |Y|^{1-p}$ and so
\begin{equation} \label{ups}
\log |\lambda| < p \, \log 2 - (p-1) \, \log |Y|.
\end{equation}
It follows from (\ref{hoop}) that $|\lambda| < 1/3$ whereby there exists $b_0 \in \mathbb{Z}$ with $|b_0| \leq (p+1) B$ such that, if we define a corresponding linear forms in logarithms
$$
\Lambda = \left| b_0 \log (-1) + b_1 \log \left( \frac{\epsilon_{3,1}}{\epsilon_{2,1}} \right) + \cdots + 
b_p \log \left( \frac{\epsilon_{3,p}}{\epsilon_{2,p}} \right) - \log \left( \zeta + 1 \right) \right|,
$$
we have $\Lambda \leq 2 | \lambda |$. We will apply Theorem \ref{Matveev}, with 
$$
n=p+2, \; D=2p(p-1) \; \mbox{ and } \;  \chi = 2.
$$
Using that
$$
h(a/b) \leq h(a)+h(b) \; \mbox{ and } \; h(a+b) \leq \log 2 + h(a) + h(b),
$$
for algebraic $a$ and $b$,
and the fact that, for each $j$, 
$$
h \left( \epsilon_{1,j} \right) = h \left( \epsilon_{2,j} \right) = h \left( \epsilon_{3,j} \right),
$$
we thus have
$$
h \left( \frac{\epsilon_{3,j}}{\epsilon_{2,j}} \right) \leq h \left( \epsilon_{3,j} \right)  +
h \left( \epsilon_{2,j} \right) \leq 2 \, h \left( \epsilon_{1,j} \right)
\; \mbox{ and } \; h (\zeta + 1) \leq \log 2.
$$
We may thus, in the notation of Theorem \ref{Matveev}, take $A_i$ such that
\begin{equation} \label{lunk}
\Omega = A_1 A_2 \cdots A_n = \pi \, \log (2) \, \left( 2 p (p-1) \right)^{p+1} \, 2^p \, \prod_{1 \leq j \leq p} h ( \epsilon_{1,j})
\end{equation}
and, via (\ref{why}),
$$
B < (p+1)  p^{2p} \, \log |Y|.
$$
We conclude therefore that 
$$
\log |\lambda| \geq \log \Lambda - \log 2 > - \log 2 - C(p+2) C_0 W_0 (2 p (p-1))^2 A_1 A_2 \cdots A_n,
$$
where, after a short computation, using the assumption that $p \geq 19$, we have
$$
C(p+2) <  (4p)^{p+4}, \; \; C_0 <  7 p \; \; \mbox{ and } \; \; W_0 < \log \log |Y| + 3 p \log p.
$$
Applying (\ref{hoop}), (\ref{lunk}) and the fact that $19 \leq p \leq 941$, we thus have 
$$
\log |\lambda| > - 2^{4p+15} \, p^{2p+8} \, (p-1)^{p+3} \, 
\left( \log \log |Y| + 3 p \log p \right) \, \prod_{1 \leq j \leq p} h ( \epsilon_{1,j}).
$$
Inequality (\ref{ilk}) thus implies that
\begin{equation} \label{ricky}
\log |\lambda| > - 2^{2p+16} \, p^{2p+11} \,  \log \left( p^{3p} \, \log |Y| \right) \,
 (p!)^2 \, R_K.
\end{equation}

We apply Lemma 9.1 of \cite{BMS} to bound the regulator $R_K$. If we suppose that we have $D_K \leq L$, then
$$
R_K < \min \{ f_K (L, 2 - t/1000) \; : \; t=0, 1, \ldots, 999 \},
$$
where
$$
f_K (L,s) = 2^{-1} \left( 2^{1-p} \pi^{-p} \sqrt{L} \right)^s \, (\Gamma (s/2))^2 \, (\Gamma (s))^{p-1} \,  s^{2p+1} \, (s-1)^{1-2p}.
$$
Since $D_K = 2^{4p-1} p^{2p}$, a short Maple computation reveals that, for $19 \leq p \leq 941$, we always have
\begin{equation} \label{trouble}
\log R_K < 10458,
\end{equation}
where the largest value of $\min \{ f_K (L, 2 - t/1000) \; : \; t=0, 1, \ldots, 999 \}$ encountered corresponds to $p=941$ and $t=743$.

Combining  (\ref{ups}), (\ref{ricky}) and (\ref{trouble}), we thus have
\begin{equation} \label{final}
\begin{array}{c}
\log |Y| < 2^{2p+17} \, p^{2p+10} \, \log \left( p^{3p} \log |Y| \right)  \,  (p!)^2 \exp(10458).
\end{array}
\end{equation}
Since $p \leq 941$,  in all cases we may conclude that $\log |Y| < 10^{15528}$, contradicting (\ref{hoop}). We may thus conclude that $(1,\pm1)$ are the only integer solutions to~\eqref{eqn:main} for the primes $p$ under consideration.

%----------------------------------------------------------------
%\subsection{The case $p=3$ (and other small primes)}
\subsection{Handling the values of $3 \leq p \leq 17$}

%----------------------------------------------------------------

%Computer algebra packages are currently somewhat ill-equipped to treat Thue equations with irrational coefficients, as in (\ref{gut}). We can get around this difficulty by reducing the problem of treating equation ~\eqref{eqn:main} for a fixed odd prime $p$ to that of solving (a finite collection of) Thue equations over $\mathbb{Z}$.
To finish the proof of Theorem~\ref{main}, it remains to solve~\eqref{eqn:main} for primes $p$ with $3 \leq p \leq 17$, which we shall carry out in this subsection.
Our strategy is to reduce the problem of treating equation ~\eqref{eqn:main} for a fixed odd prime $p$ to that of solving (a finite collection of) Thue equations over $\mathbb{Z}$.

As before, we write $\epsilon =1+\sqrt{2}$ for the  fundamental unit of $\Z[\sqrt{2}]$ (the ring of integers of $\Q(\sqrt{2})$).
Consider a solution $(x,y,p)$ to~\eqref{eqn:main} and write 
\begin{equation}\label{eqn:zDef}
z:=(x^p+3)/4,
\end{equation}
so that the integers $y$ and $z$  satisfy
$$
y^2-2z^2=-1,
$$
hence, as before (see \eqref{eqn:rel1}),
$$
y+z\sqrt{2}=s \epsilon^{k}
$$
for some odd integer $k$ and $s \in \{-1,1\}$. Replacing $y$ by $-y$ leads to replacing $k$ by $-k$ (see the proof of Lemma \ref{lem:choice}), so after possibly changing the sign of $y$ we have
\[k \equiv s \pmod{4}.\]
Writing
$$
a+b \sqrt{2}=\epsilon^{\frac{k-1}{2}},
$$
for $a, b \in \mathbb{Z}$,
we have
\begin{equation}\label{eqn:abPell}
a^2-2b^2=(-1)^{\frac{k-1}{2}}=s
\end{equation}
and
\begin{align*}
y+z\sqrt{2} & = s \epsilon (a+b \sqrt{2})^2 \\
 & = s (a^2+4ab+2b^2) +s (a^2+2ab+2b^2)\sqrt{2}.
\end{align*}
Using this parametrization for $z$ together with ~\eqref{eqn:zDef} and ~\eqref{eqn:abPell}, we have
\begin{align*}
sx^p & = 4sz-3s \\
& = 4 (a^2+2ab+2b^2)-3(a^2-2b^2)\\
& = a^2+8ab+14b^2\\
& = (a+4b)^2-2b^2.
\end{align*}
Since $a$ and $p$ are odd,  this implies  that
\begin{equation}\label{eqn:abParam}
(a+4b)+b\sqrt{2}=\epsilon^t(u+v\sqrt{2})^p
\end{equation}
for certain $u,v,t \in \Z$ with $|t| \leq (p-1)/2$.
Let us define binary forms over $\Z$ via
$$
F_{p,t}(U,V)+G_{p,t}(U,V) \sqrt{2}=\epsilon^t(U+V\sqrt{2})^p
$$
and
\[H_{p,t}(U,V)=(F_{p,t}(U,V)-4G_{p,t}(U,V))^2-2G_{p,t}(U,V)^2.\]
Then ~\eqref{eqn:abPell} and ~\eqref{eqn:abParam}  lead to
\[H_{p,t}(u,v)=s,\]
where $|u^2-2v^2| = |x|$.
We have now reduced the solution of~\eqref{eqn:main} for a fixed odd prime $p$ to the solution of the $2p$ Thue equations of degree $2p$
\begin{equation}\label{eqn:ThueSystemLarge}
H_{p,t}(U,V)=s, \quad U,V \in \Z, \quad t=0,\pm 1, \ldots, \pm \frac{p-1}{2}, \quad s \in \{-1,1\}.
\end{equation}
Taking into account information from the modular method, namely Lemma \ref{lem:sign}, we can restrict to $s=1$ when $p\geq 5$. If $p=3$, then the three
 Thue equations in \eqref{eqn:ThueSystemLarge} corresponding to $s=-1$ can readily be seen to have no solutions; for $t=-1,0,1$ there are no solutions modulo $7,9,5$ respectively.
This means that for an odd prime $p$ we only have to solve the $p$ Thue equations 
\begin{equation}\label{eqn:ThueSystem}
H_{p,t}(U,V)=1, \quad U,V \in \Z, \quad t=0,\pm 1, \ldots, \pm \frac{p-1}{2}.
\end{equation}

%%%%%%%%%%%%%%%%%%%%%%%%%%%

If we define $h_{p,t}(x)=H_{p,t}(x,1)$, let $\gamma$ be a root of $h_{p,t}(x,1)=0$ and set $K = \mathbb{Q}(\gamma)$, then we have that the discriminant of $K$ is $2^{p(6p-2)-1} p^{2p}$ and that $K$ has precisely $2$ real embeddings and hence $p$ fundamental units.

%We will appeal to a theorem of Bugeaud and Gy\H{o}ry  (part of Theorem 3 of \cite{BG}) :
%\begin{thm} \label{buggy} (Bugeaud and Gy\H{o}ry)
%Let $F(x,y)$ be an irreducible binary form with integer coefficients and degree $n \geq 3$ and let $m$ be a nonzero integer. If there exist integers $x$ and $y$ such that we have
%$$
%F(x,y)=m,
%$$
%then 
%$$
%\max \{ |x|, |y| \} < \exp \left\{ 3^{3 (n+9)} n^{18(n+1)} H^{2n-2} (\log H)^{2n-1} \log B \right\},
%$$
%where $H$ is the maximum modulus of the coefficients of $F$ and $B = \max \{|m|, e \}$.
%\end{thm}
%Applying this to (\ref{eqn:ThueSystem}) implies that we have
%$$
%\max \{ |U|, |V| \} < \exp \left\{  3^{3 (2p+9)} (2p)^{18(2p+1)} H^{4p-2} (\log H)^{4p-1} \right\}.
%$$
We note that $H_{p,0}$ is monic in $U$, so the equation $H_{p,0}(U,V)=1$ always has the solutions $U=\pm1, V=0$. If these are the only solutions to $H_{p,0}(U,V)=1$ and the other $p-1$ Thue equations in system~\eqref{eqn:ThueSystem} have no solutions, then it readily follows that $(x,y)=(1,\pm 1)$ are the only integer solution to~\eqref{eqn:main} for our fixed odd prime $p$.

We restrict our attention to $3 \leq p \leq 17$.
First of all, the Thue equation solver in {\tt PARI/GP} \cite{PARI2} can (unconditionally) solve 
$$
H_{p,0}(U,V)=1
$$
for all these primes $p$ rather quickly. The upshot of such a computation is that, apart from $(U,V)=(\pm 1, 0)$, there are no further solutions. Many of the other Thue equations in~\eqref{eqn:ThueSystem} can also be solved using {\tt PARI/GP}. However, for some of the larger values of $p$ and some values of $t$ the computations appear to take a huge amount of time (and possibly memory). Luckily many of the Thue equations involved have local obstructions to solutions (which is something not automatically checked by {\tt PARI/GP} when trying to solve the equations). In fact, the only pairs $(p,t)$ with $3 \leq p \leq 17$ and $1 \leq |t|\leq (p-1)/2$ for which we cannot find a local obstruction for $H_{p,t}(U,V)=1$ are given by
\[(p,t) \in \{(5,1), (13,1), (17,1)\}.\]
For these pairs, we can in fact use the Thue equation solver in {\tt PARI/GP} to show within reasonable time that $H_{p,t}(U,V)=1$ has no integer solutions.
For sake of completeness, for all $(p,t)$ under consideration for which we found a local obstruction, we provide in Table~\ref{table:obstructions} at least one modulus $m$ such that  the congruence $H_{p,t}(U,V)\equiv 1 \pmod{m}$ has no solutions.

 \begin{table}[h!]
 \caption{$(p,t,m)$ for which $H_{p,t}(U,V)\equiv 1 \pmod{m}$ has no solutions.}
\begin{tabular}{c|c|c}
$p$ & $t$ & $m$  \\
\hline
$3$ & $-1$ & $7$ \\
$3$ & $1$ & $3^2$ \\
$5$ & $-2,-1,2$ & $5^2$\\
$7$ & $-3,-2,-1,2,3$ & $7^2$\\
$7$ & $1,3$ & $13$\\
$11$ & $-4,-1,1,3,5$ & $11^2$ \\
$11$ & $-5,-4,-3,-2,-1,2,3,4,5$ & $23$\\
$13$ & $-6,-4,-3,4,5$ & $79$ \\
$13$ & $-5,-2,-1,2,3,6$ & $313$\\
$17$ & $-8, \ldots, -1, 2, \ldots, 8$ &$103$\\
\end{tabular}
\label{table:obstructions}
\end{table}

\begin{rem}
It is not strictly necessary to invoke the modular method to deal with the odd primes $p\leq 17$. Most of the Thue equations in~\eqref{eqn:ThueSystemLarge} with $s=-1$ can easily be solved, just like in the case when $p=3$. In fact, for all odd primes $p \leq 17 $ and all integers $t$ with $|t| \leq (p-1)/2$ we can find a local obstruction for
\begin{equation}\label{eqn:ThueMinus}
H_{p,t}(U,V)=-1
\end{equation}
(although some of the moduli involved are much larger than those in Table~\ref{table:obstructions}), except when $(p,t) \in \{(13,-4),(13,5)\}$. In the latter case the Thue equation solver in  {\tt PARI/GP} can show again that~\eqref{eqn:ThueMinus} has no solutions. For $(p,t)=(13,-4)$ we were not immediately able to solve~\eqref{eqn:ThueMinus}. We did not pursue this however, but note instead that it seems possible to use the method of Chabauty-Coleman to deal completely with~\eqref{eqn:main} for $p=13$ without using the modular method; see below.
 \end{rem}

\subsubsection{Alternative approach to handling small $p$}

To solve~\eqref{eqn:main} for a fixed odd prime $p$, without first reducing to Thue equations, we may note that this equation defines a hyperelliptic curve $C_p$ of genus $p-1$. Hence finding all integer or rational points on $C_p$ would suffice.
For convenience, consider the model for $C_p$ in weighted projective space given by
\[C_p \; :  \; 8y^2=x^{2p}+6x^p z^p+z^{2p}.\]
We want to show that $C_p(\Q)=\{(1:\pm 1: 1 ) \}$.
We have the two (non-hyperelliptic) involutions on $C_p$
\[\iota^{\pm}_p: (x:y:z) \mapsto (z:\pm y:x),\]
and their corresponding quotients
\[D^{\pm}_p:=C_p/\iota^{\pm}_p.\]
A priori, for $p>3$, it suffices to determine either $D^+_p(\Q)$ or $D^-_p(\Q)$ since $D^+_p$ and $D^-_p$ are (hyperelliptic) curves of genus $(p-1)/2>1$. For $p=3$ it would suffice of course to find that at least one of $D^+_p(\Q)$ or $D^-_p(\Q)$ is finite.

To calculate explicit models for $D^+_p$ and $D^-_p$, we find the binary forms $F^{\pm}_p$ over $\Z$ of degree $p$ such that
\[F^{\pm}_p(xz,(x \pm z)^2)=x^{2p}+6 x^p z^p+z^{2p}.\]
Introducing the variables 
\[X^{\pm}:=xz,\quad  Y^{\pm}:=(x \pm z) y, \quad Z^{\pm}:=(x \pm z)^2,\]
we see that models for $D^+_p$ and $D^-_p$ in weighted projective space are given by
\[D^{\pm}_p\; : \;  8{Y^{\pm}}^2=Z^{\pm} F_p^{\pm}(X^{\pm},Z^{\pm}).\]
Furthermore, the rational points $(1:\pm 1:1)$ on $C_p$ map under $\iota^+_p$ to $(1:\pm 2:4)$ on $D^+_p$ and under $\iota^-_p$ to $(1:0:0)$ on $D^-_p$. We note that the point at infinity $(1:0:0)$ is also a rational point on $D^+_p$. It readily follows that if we can show that
\begin{equation}\label{eqn:Dplus}
D^+_p(\Q)=\{(1:\pm 2:4), (1:0:0)\}
\end{equation}
 or
 \begin{equation}\label{eqn:Dminus}
 D^-_p(\Q)=\{(1:0:0)\},
 \end{equation}
 then $C_p(\Q)=\{(1,\pm 1, 1)\}$ and consequently $(1,\pm 1)$ are the only integer solutions to~\eqref{eqn:main} for the odd prime $p$ under consideration.
 
Using {\tt Magma}'s implementation of $2$-descent on hyperelliptic Jacobians we obtain upper bounds for the $\Q$-ranks of $\Jac(D^{\pm}_p)$ for several primes $p$; see Table~\ref{table:ranks}.
 \begin{table}[h!]
 \caption{Upper bounds for ranks of the Jacobian of $D^+_p$ and $D^-_p$.}
\begin{tabular}{c|c|c|c|c|c|c}
$p$ & 3 & 5 & 7 & 11 & 13 & 17 \\
\hline
Upper bound for $\operatorname{rank}(\Jac(D^+_p)(\Q))$ & 1 & 1 & 1 & 2 & 1 & 2 \\
\hline
Upper bound for $\rank(\Jac(D^-_p)(\Q))$ & 0 & 0 & 1 & 1 & 1 & 2 \\
\end{tabular}
\label{table:ranks}
\end{table}
We see that the ranks of $\Jac(D^-_3)(\Q)$ and $\Jac(D^-_5)(\Q)$ are both $0$. It is also easy to check that they both have trivial torsion. As a quick corollary we now obtain that~\eqref{eqn:Dminus} holds for $p=3$ and $p=5$. In fact (the Jacobian of) $D^-_3$ is the elliptic curve with Cremona Reference $1728j1$; for $D^+_3$, the reference is $1728e1$.
We can also easily check, for all primes $p$ in Table~\eqref{table:ranks}, that (two of) the three obvious rational points on $D^+_p$ give rise to a non-torsion element on the Jacobian.  We conclude that the ranks of $\Jac(D^+_7)(\Q)$ and $\Jac(D^+_{13})(\Q)$ are both $1$ and that for both of them we have an explicit generator for a subgroup of finite index. This means that it should be possible to use the method of Chabauty-Coleman to check that~\eqref{eqn:Dplus} holds for $p=7$ and $p=13$. Determining $C_p(\Q)$ for $p=11$ or $p \geq 17$ seems much harder.
%, if possible at all with current technology.

%---------------------------------------------------------------------------------------------------
\section{Frey Curves for shifted powers in more general Lucas sequences}
%----------------------------------------------------------------------------------------------------

 In this section, we will indicate how our preceding arguments fit into a more general framework. 
Let $K$ be a real quadratic number field, $\mathcal{O}_K$ its ring of integers and $\epsilon \in \mathcal{O}_K$ a fundamental unit in $K$, with conjugate $\overline{\epsilon}$. Define the Lucas sequences, of the first and second kinds, respectively, 
$$
U_k=\frac{\epsilon^k-\left(\overline{\epsilon}\right)^k}{\epsilon-\overline{\epsilon}} \; \mbox{ and } \;
V_k =  \epsilon^k+\left(\overline{\epsilon}\right)^k, \; \mbox{ for } \;  k \in \Z.
$$
Let $a,c \in \Q$ with $a\not=0$, and consider the problem of determining the shifted powers $ax^n+c$ in one of these sequences, i.e. determining all integers $k, x$ and $n$ with $n \geq 2$ such that we have
\begin{equation}\label{eqn:ShiftedPowers}
U_k=ax^n+c
\end{equation}
or
\begin{equation}\label{eqn:ShiftedPowers2}
V_k=ax^n+c.
\end{equation}
If e.g. $\epsilon=(1+\sqrt{5})/2$, this amounts to determining  shifted powers in the Fibonacci ($U_k$) or Lucas ($V_k$) sequences. If $\epsilon=1+\sqrt{2}$, $(a,c)=(\pm 1/4, \pm 3/4)$ and $k$ is odd, equation (\ref{eqn:ShiftedPowers}) corresponds to  the main problem of this paper.

We will show that the arguments of this paper can potentially resolve such problems corresponding to either

\begin{itemize}
\item equation (\ref{eqn:ShiftedPowers}) with $k$ odd  and $\Norm(\epsilon)=-1$, or
\item equation (\ref{eqn:ShiftedPowers2}) with either $k$ even or  $\Norm(\epsilon)=1$.
\end{itemize}
We proceed from the observation that
\begin{equation} \label{crux}
\epsilon^k + \epsilon^{-k} \pm 2 =
\left\{
\begin{array}{ll}
\left( \epsilon^{k/2} \pm \epsilon^{-k/2} \right)^2 & \mbox{ if $k$ is even, } \\
\epsilon \left(  \epsilon^{\frac{k-1}{2}} \pm \epsilon^{\frac{-k-1}{2}} \right)^2 & \mbox{ if $k$ is odd. } \\
\end{array}
\right.
\end{equation}
If $k$ is odd  and $\Norm(\epsilon)=-1$, we have $\left(\overline{\epsilon}\right)^k=-\epsilon^{-k}$
whereby
$$
\left( \epsilon+ \epsilon^{-1} \right) U_k \pm 2 = \epsilon \left(  \epsilon^{\frac{k-1}{2}} \pm \epsilon^{\frac{-k-1}{2}} \right)^2.
$$
It follows from equation (\ref{eqn:ShiftedPowers}) that
\begin{equation} \label{tom2}
(\epsilon+\epsilon^{-1})ax^n+((\epsilon+\epsilon^{-1} )c\pm 2)=\epsilon \gamma_{k,\pm}^2,
\end{equation}
where
$$
\gamma_{k,\pm}:=\epsilon^{\frac{k-1}{2}}\pm \epsilon^{\frac{-k-1}{2}} \in \mathcal{O}_K,
$$
Similarly, if either $k$ is even, or $\Norm(\epsilon)=1$, we have
$$
V_k = \epsilon^k + \epsilon^{-k},
$$
and hence, from (\ref{crux}) and assuming that we have (\ref{eqn:ShiftedPowers2}), 
\begin{equation} \label{tom3}
ax^n+(c\pm 2)=
\left\{
\begin{array}{ll}
\tau_{k,\pm}^2 & \mbox{ if $k$ is even, } \\
\epsilon \gamma_{k,\pm}^2 & \mbox{ if $k$ is odd, } \\
\end{array}
\right.
\end{equation}
where
$$
\tau_{k,\pm}:=\epsilon^{k/2}\pm \epsilon^{-k/2} \in \mathcal{O}_K.
$$

To either three term relation (\ref{tom2}) or (\ref{tom3}), we can actually associate several Frey curves -- the most obvious are those corresponding  to the generalized Fermat equations of signature $(n,n,2)$ or signature $(n,3,2)$. Both these Frey curves are defined over the totally real number field $K$ and therefore amenable to a Hilbert modular approach. In this setting, we can in fact generalize our choices of $a,c \in \Q$ to $a,c \in K$ and $x \in \Z$ to $x \in \mathcal{O}_K$.

\begin{rem}
If we also want to study~\eqref{eqn:ShiftedPowers} for $k$ even or $\Norm(\epsilon)=1$, or equation (\ref{eqn:ShiftedPowers2}) when $k$ is odd  and $\Norm(\epsilon)=-1$,  we can write down  similar three term relations to those above, but now over $K(\sqrt{-1})$. The resulting Frey curves will therefore not {\it a priori} be defined over a totally real number field.
\end{rem}

%---------------------------------------------------------
\section{Other applications}
%---------------------------------------------------------

In \cite{BLMS}, equation (\ref{cur}) is solved in the case where $\{ u_k \} = \{ F_k \}$, the Fibonacci numbers and $c = \pm 1$. In this situation, the problem actually reduces to that of determining (almost) perfect powers in the Fibonacci and Lucas sequences (a program that is carried out in \cite{BMS}; see also \cite{BLMS2}), through a series of identities akin to
$$
F_{4k}+1 = F_{2k-1} L_{2k+1}.
$$
Here, $\{ L_k \}$ denotes the Lucas numbers, the companion sequence to the Fibonacci numbers.
This reduction in case $c = \pm 1$ depends crucially upon the fact that
$$
F_{-1}=F_1=F_2=1 \; \mbox{ and } \; F_{-2}=-1,
$$
and does not apparently extend to permit solution of the equation $F_k=ax^n+c$ for even a single fixed pair $(a,c)$ with $|c|>1$ (the case $c=-2$ is given as an open problem in \cite{BLMS}).

As explained in the preceding section, the methods of this paper potentially permit solution of $F_k=ax^n+c$ for any given pair $(a,c)$ and  {\it odd} values of the index $k$. For certain choices of $c$, however, the elementary arguments of \cite{BLMS} allow one to handle the remaining even index terms. By way of example, let us suppose that $c = F_{2j}$, for some integer $j$, so that
$$
|c| \in \left\{ 1, 3, 8, 21, 55, 144, \ldots  \right\}
$$
whereby we are considering equations of the shape
\begin{equation} \label{sharknado2}
F_k \pm F_{2j} = a x^n,
\end{equation}
for given fixed $a$ and $j$ (here, $k, x$ and $n$ are variables).
Then, appealing to the identity
$$
F_i L_j = F_{i+j} + (-1)^j  F_{i-j},
$$
we have that
$$
F_{4k} = F_{2k+j} L_{2k-j} + (-1)^{j+1} F_{2j}
$$
and
$$
F_{4k} = F_{2k-j} L_{2k+j} + (-1)^{j+1} F_{-2j} = F_{2k-j} L_{2k+j} + (-1)^{j} F_{2j}.
$$
Assuming that we have $F_{4k}=a x^n+c$, if also $j$ is even, say, it follows that
$$
a x^n = F_{2k-j} L_{2k+j},
$$
whilst if $j$ is odd,
$$
a x^n = F_{2k+j} L_{2k-j}.
$$
Similarly, the identities
$$
F_{4k+2} = F_{2k+j+1} L_{2k-j+1} + (-1)^{j} F_{2j}
$$
and
$$
F_{4k+2} = F_{2k-j+1} L_{2k+j+1} + (-1)^{j+1} F_{2j},
$$
with $F_{4k+2}=a x^n+c$, imply that
$$
a x^n = F_{2k+j+1} L_{2k-j+1} \; \mbox{ or } \; a x^n = F_{2k-j+1} L_{2k+j+1}.
$$
In all cases, we are able to descend to a problem of (almost) perfect powers in the Fibonacci or Lucas sequences, one that has proven to be computationally tractable (whereby the same is potentially true for equation (\ref{sharknado2})).

%%%%%%%%%%%%%%%%%%%%%%%%%%%

\end{document}